\documentclass[a4paper]{elsarticle}
\usepackage{multirow}
\usepackage[fleqn]{amsmath}
\usepackage{graphicx}
\usepackage{caption}
\usepackage{dsfont}
\usepackage{amsfonts}
\usepackage{subfigure}
\usepackage{lineno,hyperref}
\modulolinenumbers[5]
\usepackage{CJK}
\usepackage{indentfirst}
\usepackage{amsmath}
\usepackage{color}

\usepackage{geometry}
\journal{Computers and Mathematics with Applications}

\begin{document}

\begin{frontmatter}


\title{A Galerkin finite element method for time-fractional stochastic heat equation}

\author[mymainaddress]{Guang-an Zou\corref{mycorrespondingauthor}}
\address[mymainaddress]{School of Mathematics and Statistics, Henan University, Kaifeng 475004, P. R. China}
\cortext[mycorrespondingauthor]{Corresponding author}
\ead{zouguangan00@163.com,zouguangan@henu.edu.cn}

\begin{abstract}

In this study, a Galerkin finite element method is presented for time-fractional stochastic heat equation driven by multiplicative noise, which arises from the consideration of heat transport in porous media with thermal memory with random effects. The spatial and temporal regularity properties of mild solution to the given problem under certain sufficient conditions are obtained. Numerical techniques are developed by the standard Galerkin finite element method in spatial direction, and Gorenflo-Mainardi-Moretti-Paradisi scheme is applied in temporal direction.  The convergence error estimates for both semi-discrete and fully discrete schemes are established. Finally, numerical example is provided to verify the theoretical results.

\end{abstract}


\end{frontmatter}


\section{Introduction}

Over the last few decades, fractional calculus (i.e., fractional integrals and fractional derivatives) have attracted considerable interests primarily due to their potential applications in various fields of science and engineering [1,2,40,41]. As we know, fractional differential equations are highly effective mathematical tools to describe the complex behaviors and phenomena of memory processes [1,10,15,35], it also can effectively characterize the ubiquitous power-law phenomena [36]. Many theoretical analysis and numerical methods are developed for fractional differential equations, see the literatures [9,10,12,16-22,25-35] and the references therein. On the other hand, stochastic perturbations are coming from many natural sources in the practically physical system, they can not be ignored and the presence of noises might give rise to some statistical features and important phenomena, then the stochastic differential equations are produced, which are more realistic mathematical model of the real-world situations [8]. Recently, some related works about the theoretical analysis of fractional stochastic differential equations have been intensively investigated in the literatures [3-7,11,23,24,43]. However, it seems that there are less literatures related to numerical approximation of stochastic partial differential equations with fractional derivatives.

In this article, we consider the following time-fractional stochastic heat equation perturbed by multiplicative noise on the finite time interval $[0,T]$:
\begin{align*}
 ^{C}D_{t}^{\alpha}u(x,t)=\Delta u(x,t)+\sigma(u(x,t))\frac{d W(t)}{d t},x\in D,t\in(0,T], \tag{1.1}
\end{align*}
subject to the initial condition:
\begin{align*}
 u(x,0)=u_{0}(x),x\in D,t=0, \tag{1.2}
\end{align*}
and the boundary conditions
\begin{align*}
 u(x,t)=0,x\in \partial D,t>0, \tag{1.3}
\end{align*}
where $u(x,t)$ is a random function, $D$ is a bounded subset of $\mathbb{R}^{d}$ (where $d=1,2,3$) with the boundary $\partial D$; the operator $\Delta$ stands for the Laplacian acting on $\mathbb{R}^{d}$. The coefficient $\sigma$ is real-valued continuous function; denote by $W(t)$ the Wiener process given on a filtered probability space $(\Omega,\mathcal{F},\mathbf{P},\{\mathcal{F}_{t}\}_{t\geq 0})$. Here, we define the Caputo derivative of order $\alpha$ as (see Ref.[1])
\begin{align*}
 ^{C}D_{t}^{\alpha}u(x,t)=\begin{cases}
\frac{1}{\Gamma(n-\alpha)}\int_{0}^{t}(t-s)^{n-\alpha-1}\frac{\partial^{n}u(x,s)}{\partial s^{n}}ds,n-1< \alpha < n,\\
\frac{\partial^{n}u(x,t)}{\partial t^{n}}, \qquad \qquad \qquad \qquad \qquad \quad \alpha = n,\\
\end{cases} \tag{1.4}
\end{align*}
in which we denote $n=[\alpha]+1$ wherein $[\alpha]$ representing the integer part of $\alpha$, and $\Gamma(\cdot)$ stands for Gamma function $\Gamma(\alpha)=\int_{0}^{\infty}t^{\alpha-1}e^{-t}dt$.

Noticing that the deterministic time-fractional heat equation (where $\sigma=0$ in Eq.(1.1)) describe the diverse anomalous diffusive processes in complex media with different $\alpha$, for example, subdiffusion ($0<\alpha<1$), normal diffusion ($\alpha=1$), superdiffusion ($1<\alpha<2$), and ballistic diffusion ($\alpha=2$)(see [2]). In this paper, we focus on the fractional case $0<\alpha<1$, then the operator $^{C}D_{t}^{\alpha}$ can be written as
\begin{align*}
 ^{C}D_{t}^{\alpha}u(x,t)=\frac{1}{\Gamma(1-\alpha)}\int_{0}^{t}\frac{\partial u(x,s)}{\partial s}\frac{ds}{(t-s)^{\alpha}}, 0<\alpha<1. \tag{1.5}
\end{align*}

Note that the initial-boundary value problem (1.1)-(1.3) is a special case of the fractional stochastic partial differential equations (SPDEs) discussed in [3-7], which can be used to model the random effects on transport of particles in medium with thermal memory. For details, Chen et al. [3] introduced a class of SPDEs with time-fractional derivatives and proved
the existence and uniqueness of solutions to these equations. Mijena and Nane [4] proved the existence and uniqueness of mild solution to non-linear space-time fractional SPDEs, and they also investigated the bounds for the intermittency fronts solution of these equations [5]. Foondun and Nane [6] studied the asymptotic properties of space-time fractional SPDEs. Chen et al. [7] proved the existence and uniqueness of solution to space-time fractional SPDEs in Gaussian noisy environment. Zou et al. [24] studied the existence and regularity of mild solution to this type of fractional stochastic evolution equation, they also used a semi-discrete finite element method for solving fractional stochastic diffusion-wave equations [42]. Li et al. [36] developed a Galerkin finite element approximations for stochastic space-time fractional wave equations. However, to the best of our knowledge, numerical methods for these kinds of fractional stochastic subdiffusion problem are yet to be investigated. Motivated by this facts, our goal of this paper is to develop a Galerkin finite element method for time-fractional SPDEs and demonstrate the application of this method with the aid of example.

The remaining of this paper is organized as follows. In Section 2, some notations and preliminaries are recalled, and we prove the spatial and temporal regularity properties of mild solution to time-fractional stochastic heat equation. In Section 3, we propose the semi-discrete and fully discrete finite element methods for this time-fractional stochastic evolution equation, and the strong convergence error estimates for both semi-discrete and fully discrete schemes in $L_{2}$-norm are established. Numerical example is presented in Section 4. Conclusions and discussions are given in the final section.

\section{Notations and preliminaries}

Denote by $H=L_{2}(D)$ a real separable Hilbert space with inner product $(\cdot,\cdot)$ and norm $\|\cdot\|$. We assume that $L_{2}(\Omega,H)$ is Hilbert space of $H$-valued random variables equipped with the inner product $\mathbf{E}(\cdot,\cdot)$ and norm $\mathbf{E}\|\cdot\|$, given by
\begin{align*}
L_{2}(\Omega,H)=\{v:\mathbf{E}\|v\|^{2}=\int_{\Omega}\|v(\omega)\|^{2}d\mathbf{P}(\omega)<\infty,\omega\in \Omega\},
\end{align*}
where $\mathbf{E}$ denote the expectation.

To be more specific, let $W(t)$ be a Wiener process with the linear bounded covariance operator $Q$. Moreover, there exist the eigenvalues $\lambda_{n}$ of $Q$ with corresponding eigenfunctions $e_{n}$ such that $Q e_{n}=\lambda_{n}e_{n}$ and $Tr(Q)=\sum\limits_{n=1}^{\infty}\lambda_{n}<\infty$. Then we have the following representation of $W(t)$:
\begin{align*}
W(t)=\sum\limits_{n=1}^{\infty}\lambda^{1/2}_{n}\beta_{n}(t)e_{n},
\end{align*}
where the sequence $\{\beta_{n}\}_{n\geq1}$ is real-valued standard Brownian motions.

Let us denote by $L_{2}^{0}=L_{2}(Q^{1/2}(H),H)$ the space of linear Hillbert-Schmidt operators from a Hilbert space $Q^{1/2}(H)$ to a Hilbert space $H$, which is defined by
\begin{align*}
L_{2}^{0}=\{\phi\in L(H):\sum\limits_{n=1}^{\infty}\|\phi Q^{1/2}e_{n})\|^{2}<\infty\},
\end{align*}
where $L(H)$ is space of linear bounded operators from $H$ to $H$. Now we give the following important result for the stochastic integral.

\textbf{Lemma 2.1.}([8]) Let $v:[0,T]\times\Omega\rightarrow L_{2}^{0}$ be a strongly measurable mapping such that $\int_{0}^{T}\mathbf{E}\|v(s)\|_{L_{2}^{0}}^{2}ds<\infty$. Then the following It\^{o}-isometry holds :
\begin{align*}
\mathbf{E}\|\int_{0}^{t}v(s)dW(s)\|^{2}=\int_{0}^{t}\mathbf{E}\|v(s)\|_{L_{2}^{0}}^{2}ds  \tag{2.1}
\end{align*}
for all $ 0\leq t\leq T$.

For the sake of convenience, we can rewrite the time-fractional stochastic heat problem (1.1)-(1.3) as the following abstract form:
\begin{align*}
\begin{cases}
 ^{C}D_{t}^{\alpha}u(t)=Au(t)+\sigma(u(t))\frac{d W}{d t},x\in D,t>0,\\
 u(0)=u_{0}, \end{cases}\tag{2.2}
\end{align*}
in which we denote $u(t)=u(\cdot,t)$ and $A=\Delta$ with the domain $D(A)=H^{2}(D)\bigcap H_{0}^{1}$, where $H^{2}(D)=\{v\in L_{2}(D),\frac{d^{k}v}{dx^{k}}\in L_{2}(D),k\leq 2; k\in N^{+}\}$ and $H_{0}^{1}=H_{0}^{1}(D)=\{v\in H^{1}(D), v|_{\partial D}=0 \}$.

Let $A^\mu$ be a fractional power of $A$ and the spaces $\dot{H}^{\mu}=\mathcal{D}(A^{\frac{\mu}{2}})$ be equipped with the induced norm $\|v\|_{\mu}=\|A^{\frac{\mu}{2}}v\|$. Denote by $\dot{H}^{0}=H$, it is clear that $\dot{H}^{1}=H_{0}^{1}$, $\dot{H}^{2}=H^{2}(D)\bigcap H_{0}^{1}$ with equivalent norms. The space $L_{2}(\Omega,\dot{H}^{\mu})$ equipped with norm defined by
\begin{align*}
\|v\|_{L_{2}(\Omega,\dot{H}^{\mu})}=\mathbf{E}(\|v\|_{\mu}^{2})^{\frac{1}{2}}=(\int_{\Omega}\|v(\omega)\|_{\mu}^{2}d\mathbf{P}(\omega))^{\frac{1}{2}},
\end{align*}
where $\mu\in R$ and $\omega\in \Omega$.

Throughout the paper, we impose the assumption that the measurable function $\sigma: H \rightarrow L_{2}^{0}$ satisfies the following global Lipschitz and linear growth conditions
\begin{align*}
 \|\sigma(u)-\sigma(v)\|_{L_{2}^{0}}\leq C_{\sigma}\|u-v\|, \|\sigma(v)\|_{L_{2}^{0}}\leq C_{\sigma}\|v\|,\tag{2.3}
\end{align*}
for any $u,v\in H$ and $C_{\sigma}>0$ is a constant.

Now, inspired by the definition of mild solutions to time-fractional differential equations (see Refs.[9-11,23,24]), we give the following definition of mild solution for our time-fractional stochastic heat equation.

\textbf{Definition 2.1.} An adapted process $\{u(t)\}_{t\geq0}$ is called a mild solution to (2.2), if it satisfies the following integral equation $\mathbf{P}$-a.s.
\begin{align*}
 u(t)=E_{\alpha}(t)u_{0}+\int_{0}^{t}(t-s)^{\alpha-1}E_{\alpha,\alpha}(t-s)\sigma(u(s))dW(s),  \tag{2.4}
\end{align*}
for a.s. $\omega\in\Omega$, where the generalized Mittag-Leffler operators are given by
\begin{align*}
E_{\alpha}(t)=\int_{0}^{\infty}\xi_{\alpha}(\theta)S(t^{\alpha}\theta)d\theta,
~E_{\alpha,\alpha}(t)=\int_{0}^{\infty}\alpha\theta\xi_{\alpha}(\theta)S(t^{\alpha}\theta)d\theta,
\end{align*}
which involving the Wright-type function:
\begin{align*}
\xi_{\alpha}(\theta)=\frac{1}{\alpha}\theta^{-1-\frac{1}{\alpha}}\omega_{\alpha}(\theta^{-\frac{1}{\alpha}})\geq0,~\alpha\in(0,1)
\end{align*}
connected with the following one-sided stable probability function:
\begin{align*}
\omega_{\alpha}(\theta)=\frac{1}{\pi}\sum_{n=1}^{\infty}(-1)^{n-1}\theta^{-n\alpha-1}\frac{\Gamma(n\alpha+1)}{n!}\sin(n\pi\alpha),\theta\in(0,\infty),
\end{align*}
in which $S(t)=e^{tA}$ is an analytic semigroup generated by a linear operator $A$. The proof of existence and uniqueness of mild solution to (2.2) are similar to the literatures [11,23,24].

In what follows, let us introduce and prove some important lemmas, which will be used in the
subsequent discussions.

\textbf{Lemma 2.2.} (see [10,24]) For any $\alpha\in(0,1)$ and $\nu\in (-1,\infty)$, it is not difficult to verity that
\begin{align*}
\int_{0}^{\infty}\theta^{\nu}\xi_{\alpha}(\theta)d\theta=\frac{\Gamma(1+\nu)}{\Gamma(1+\alpha\nu)}. \tag{2.5}
\end{align*}

\textbf{Lemma 2.3.} (see [13]) For any $\gamma\geq 0 $, there exists a constant $C_{\gamma}>0$ such that
\begin{align*}
\|A^{\gamma}S(t)\|\leq C_{\gamma}t^{-\gamma}, t>0.\tag{2.6}
\end{align*}

\textbf{Lemma 2.4.} For any $t>0$ and $0\leq\gamma<1$, there exist constants $C_{\alpha1}=C(\alpha,\gamma)>0$ and $C_{\alpha2}=C(\alpha,\gamma)>0$ such that
\begin{align*}
\|A^{\gamma}E_{\alpha}(t)v\|\leq C_{\alpha1} t^{-\gamma\alpha}\|v\|,~\|A^{\gamma}E_{\alpha,\alpha}(t)v\|\leq C_{\alpha2}t^{-\gamma\alpha}\|v\|.\tag{2.7}
\end{align*}

\textbf{Proof.} For $t>0$ and $0\leq\gamma<1$, using Lemma 2.2 and Lemma 2.3, we get
\begin{align*}
\|A^{\gamma}E_{\alpha}(t)v\|&\leq\int_{0}^{\infty}\xi_{\alpha}(\theta)\|A^{\gamma}S(t^{\alpha}\theta)v\|d\theta\\
&\leq\int_{0}^{\infty}C_{\gamma}t^{-\gamma\alpha}\theta^{-\gamma}\xi_{\alpha}(\theta)\|v\|d\theta\\
&=\frac{C_{\gamma}\Gamma(1-\gamma)}{\Gamma(1-\gamma\alpha)}t^{-\gamma\alpha}\|v\|, v\in L_{2}(D),
\end{align*}
and
\begin{align*}
\|A^{\gamma}E_{\alpha,\alpha}(t)v\|&\leq\int_{0}^{\infty}\alpha\theta \xi_{\alpha}(\theta)\|A^{\gamma}S(t^{\alpha}\theta)v\|d\theta\\
&\leq \int_{0}^{\infty}C_{\gamma}\alpha t^{-\gamma\alpha}\theta^{1-\gamma}\xi_{\alpha}(\theta) \|v\|d\theta\\
&=\frac{C_{\gamma}\alpha\Gamma(2-\gamma)}{\Gamma(1+\alpha(1-\gamma))}t^{-\gamma\alpha}\|v\|, v\in L_{2}(D).
\end{align*}

This completes the proof.  $\square$

\textbf{Remark 2.1.} If $\gamma=0$, it is easy to see that $E_{\alpha}(t)$ and $E_{\alpha,\alpha}(t)$ are linear and bounded operators.

\textbf{Lemma 2.5.} For any $0\leq t_{1}< t_{2}\leq T$ and $0\leq a\leq1$, we have
\begin{align*}
t_{2}^{a}-t_{1}^{a}\leq a\delta^{a-1}(t_{2}-t_{1})^{a},\tag{2.8}
\end{align*}
where $0<\delta<1$.

\textbf{Proof.} This inequality can be proved by means of Lagrange's mean value theorem, we omit it here.

\textbf{Lemma 2.6.} For any $0\leq t_{1}< t_{2}\leq T$ and for $0<\gamma<1$, there exist two constants $C_{\gamma1}=C(\alpha,\gamma)>0$ and $C_{\gamma2}=C(\alpha,\gamma)>0$ such that
\begin{align*}
\|A^{\gamma}[E_{\alpha}(t_{2})-E_{\alpha}(t_{1})]v\|\leq C_{\gamma1}(t_{2}-t_{1})^{\gamma\alpha}\|v\|,\tag{2.9}
\end{align*}
and
\begin{align*}
\|A^{\gamma}[E_{\alpha,\alpha}(t_{2})-E_{\alpha,\alpha}(t_{1})]v\|\leq C_{\gamma2}(t_{2}-t_{1})^{\gamma\alpha}\|v\|.\tag{2.10}
\end{align*}

\textbf{Proof.} For any $0\leq t_{1}<t_{2}\leq T$, we can deduce that
\begin{align*}
S(t_{2}^{\alpha}\theta)-S(t_{1}^{\alpha}\theta)=\int_{t_{1}}^{t_{2}}\frac{d S(t^{\alpha}\theta)}{dt} dt=-\int_{t_{1}}^{t_{2}} \alpha t^{\alpha-1}\theta A S(t^{\alpha}\theta)dt.\tag{2.11}
\end{align*}

For $0<\gamma<1$, making use of the expression (2.11), by Lemma 2.2, Lemma 2.3 and Lemma 2.5, we obtain
\begin{align*}
\|A^{\gamma}[E_{\alpha}(t_{2})-E_{\alpha}(t_{1})]v\|&=\|\int_{0}^{\infty}\xi_{\alpha}(\theta)A^{\gamma}[S(t_{2}^{\alpha}\theta)-S(t_{1}^{\alpha}\theta)]vd\theta\|\\
&\leq\int_{0}^{\infty}\alpha\theta \xi_{\alpha}(\theta)\int_{t_{1}}^{t_{2}} t^{\alpha-1}\| A^{\gamma+1} S(t^{\alpha}\theta)v\|dtd\theta\\
&\leq\int_{0}^{\infty}C_{\gamma}\alpha\theta^{-\gamma}\xi_{\alpha}(\theta)(\int_{t_{1}}^{t_{2}}t^{-\gamma\alpha-1}dt) \|v\|d\theta\\
&=\frac{ C_{\gamma}\Gamma(1-\gamma)}{\gamma\Gamma(1-\gamma\alpha)}(t_{1}^{-\gamma\alpha}-t_{2}^{-\gamma\alpha})\|v\|\\
&\leq \frac{C_{\gamma}^{\dagger}\Gamma(1-\gamma)}{\gamma T_{0}^{2\gamma\alpha}\Gamma(1-\gamma\alpha)}(t_{2}-t_{1})^{\gamma\alpha}\|v\|, v\in L_{2}(D),
\end{align*}
and
\begin{align*}
\|A^{\gamma}[E_{\alpha,\alpha}(t_{2})-E_{\alpha,\alpha}(t_{1})]v\|&=\|\int_{0}^{\infty}\alpha\theta \xi_{\alpha}(\theta)A^{\gamma}[S(t_{2}^{\alpha}\theta)-S(t_{1}^{\alpha}\theta)]vd\theta\|\\
&\leq\int_{0}^{\infty}\alpha^{2}\theta^{2} \xi_{\alpha}(\theta)\int_{t_{1}}^{t_{2}} t^{\alpha-1}\| A^{\gamma+1} S(t^{\alpha}\theta)v\|dtd\theta\\
&\leq\int_{0}^{\infty}C_{\gamma}\alpha^{2}\theta^{1-\gamma}\xi_{\alpha}(\theta)(\int_{t_{1}}^{t_{2}}t^{-\gamma\alpha-1}dt) \|v\|d\theta\\
&=\frac{C_{\gamma}\alpha\Gamma(2-\gamma)}{\gamma\Gamma(1+\alpha(1-\gamma))}(t_{1}^{-\gamma\alpha}-t_{2}^{-\gamma\alpha})\|v\|\\
&\leq \frac{C_{\gamma}^{\dagger}\alpha\Gamma(2-\gamma)}{\gamma T_{0}^{2\gamma\alpha}\Gamma(1+\alpha(1-\gamma))}(t_{2}-t_{1})^{\gamma\alpha}\|v\|, v\in L_{2}(D),
\end{align*}
where $0<T_{0}< T$. This completes the proof.  $\square$

\textbf{Remark 2.2.} The operators $E_{\alpha}(t)$ and $E_{\alpha,\alpha}(t)$ are strongly continuous. Furthermore, if $\gamma=0$ in Lemma 2.6, then there exist two constants $C_{\alpha1}>0$ and $C_{\alpha2}>0$ such that
\begin{align*}
\|(E_{\alpha}(t_{2})-E_{\alpha}(t_{1}))v\|\leq C_{\alpha1}(t_{2}-t_{1})\|v\|, \tag{2.12}
\end{align*}
and
\begin{align*}
\|(E_{\alpha,\alpha}(t_{2})-E_{\alpha,\alpha}(t_{1}))v\|\leq C_{\alpha2}(t_{2}-t_{1})\|v\|. \tag{2.13}
\end{align*}

\textbf{Proof.} For any $0\leq t_{1}<t_{2}\leq T$, the same as the proof of Lemma 2.6, we conclude that
\begin{align*}
\|(E_{\alpha}(t_{2})-E_{\alpha}(t_{1}))v\|&=\|\int_{0}^{\infty}\xi_{\alpha}(\theta)(S(t_{2}^{\alpha}\theta)-S(t_{1}^{\alpha}\theta)vd\theta\|\\
&\leq\int_{0}^{\infty}\alpha\theta \xi_{\alpha}(\theta)\int_{t_{1}}^{t_{2}} t^{\alpha-1}\| AS(t^{\alpha}\theta)v\|dtd\theta\\
&\leq \int_{0}^{\infty}C_{0}\alpha\xi_{\alpha}(\theta)(\int_{t_{1}}^{t_{2}}t^{-1}dt) \|v\|d\theta\\
&=C_{0}\alpha(\ln t_{2}-\ln t_{1})\|v\|\\
&\leq \frac{C_{0}\alpha}{T_{0}}(t_{2}-t_{1})\|v\|, v\in L_{2}(D),
\end{align*}
and
\begin{align*}
\|(E_{\alpha,\alpha}(t_{2})-E_{\alpha,\alpha}(t_{1}))v\|&=\|\int_{0}^{\infty}\alpha\theta \xi_{\alpha}(\theta)(S(t_{2}^{\alpha}\theta)-S(t_{1}^{\alpha}\theta)vd\theta\|\\
&\leq\int_{0}^{\infty}\alpha^{2}\theta^{2} \xi_{\alpha}(\theta)\int_{t_{1}}^{t_{2}} t^{\alpha-1}\| A S(t^{\alpha}\theta)v\|dtd\theta\\
&\leq\int_{0}^{\infty} C_{0}\alpha^{2}\theta \xi_{\alpha}(\theta)(\int_{t_{1}}^{t_{2}}t^{-1}dt) \|v\|d\theta\\
&=\frac{C_{0}\alpha^{2}\Gamma(2)}{\Gamma(1+\alpha)}(\ln t_{2}-\ln t_{1})\|v\|\\
&\leq \frac{C_{0}\alpha^{2}\Gamma(2)}{T_{0}\Gamma(1+\alpha)}(t_{2}-t_{1})\|v\|, v\in L_{2}(D),
\end{align*}
where $0<T_{0}< T$. It is easy to see that the operators $E_{\alpha}(t)$ and $E_{\alpha,\alpha}(t)$ are strongly continuous.

Next, we shall discuss and prove the spatial and temporal regularity properties of mild solution to Eq.(2.2). To begin with, we will introduce a generalization of standard integral version of Gronwall's Lemma with weak singularities [14,44].

\textbf{ Lemma 2.7.} Let $T>0$ and $C_{1},C_{2}>0$ and let $\varphi:[0,T]\rightarrow R$ be a nonnegative and continuous function. Let $\beta>0$. If we have
\begin{align*}
\varphi(t)\leq C_{1}+C_{2}\int_{0}^{t}(t-s)^{\beta-1}\varphi(s)ds
\end{align*}
then there exists a constant $C=C(C_{2},T,\beta)$ such that
\begin{align*}
\varphi(t)\leq CC_{1}
\end{align*}
for all $t\in(0,T]$.

\textbf{Theorem 2.1.} Assume that $\sigma$ satisfies (2.3), for any $\alpha\in (\frac{1}{2},1)$ and $0\leq \nu \leq 1$. Let $u(t)$ be a mild solution to (2.2). Then there exists a constant $C>0$ such that
\begin{align*}
\sup\limits_{t\in [0,T]}\|u(t)\|_{L_{2}(\Omega,\dot{H}^{\nu})}^{2}\leq C\|u_{0}\|_{L_{2}(\Omega,H)}^{2}. \tag{2.14}
\end{align*}
\textbf{Proof.} For any $\frac{1}{2}<\alpha<1$ and $0\leq \nu \leq 1$, from (2.4), using Lemma 2.1 and Lemma 2.4, we have
\begin{align*}
\mathbf{E}\|u(t)\|_{\nu}^{2}&\leq 2\mathbf{E}\|A^{\frac{\nu}{2}}E_{\alpha}(t)u_{0}\|^{2}+2\mathbf{E}\|\int_{0}^{t}(t-s)^{\alpha-1}A^{\frac{\nu}{2}}E_{\alpha,\alpha}(t-s)\sigma(u(s))dW(s)\|^{2}\\
&\leq 2C_{\alpha1}^{2} t^{-\nu\alpha}\mathbf{E}\|u_{0}\|^{2}+2(\int_{0}^{t}\mathbf{E}\|(t-s)^{\alpha-1}A^{\frac{\nu}{2}}E_{\alpha,\alpha}(t-s)\sigma(u(s))\|_{L_{2}^{0}}^{2}ds)\\
&\leq 2C_{\alpha1}^{2} T_{0}^{-\nu\alpha}\mathbf{E}\|u_{0}\|^{2}+2C_{\alpha2}^{2}C_{\sigma}^{2}(\int_{0}^{t}(t-s)^{(\beta-1)}\mathbf{E}\|u(s)\|_{\nu}^{2}ds),
\end{align*}
where $\beta=2\alpha-1>0$ and $0<T_{0}<T$.

By means of Lemma 2.7, then the proof of Theorem 2.1 is finished.  $\square$

\textbf{Theorem 2.2.} Assume that $\sigma$ satisfies (2.3), for any $0\leq t_{1}<t_{2}\leq T$, $\frac{1}{2}<\alpha<1$ and $0\leq \nu\leq 1$. Then a mild solution $u(t)$ to (2.2) is H\"{o}lder continuous with respect to the norm $\|\cdot\|_{L_{2}(\Omega,\dot{H}^{\nu})}$ and satisfies
\begin{align*}
\|u(t_{2})-u(t_{1})\|_{L_{2}(\Omega,\dot{H}^{\nu})}^{2}\leq C(t_{2}-t_{1})^{\beta},\tag{2.15}
\end{align*}
where $\beta=\max\{\nu\alpha,2(1-\alpha),(2\alpha-1)\}>0$. 

\textbf{Proof.} For any $0\leq t_{1}<t_{2}\leq T$, from (2.4) we get
\begin{align*}
u(t_{2})-u(t_{1})&=E_{\alpha}(t_{2})u_{0}-E_{\alpha}(t_{1})u_{0}\\
&\hspace{2mm}+\int_{0}^{t_{2}}(t_{2}-s)^{\alpha-1}E_{\alpha,\alpha}(t_{2}-s)\sigma(u(s))dW(s)\\
&\hspace{2mm}-\int_{0}^{t_{1}}(t_{1}-s)^{\alpha-1}E_{\alpha,\alpha}(t_{1}-s)\sigma(u(s))dW(s)\\
&=: I_{1}+I_{2},\tag{2.16}
\end{align*}
where
\begin{align*}
I_{1}=[E_{\alpha}(t_{2})-E_{\alpha}(t_{1})]u_{0},
\end{align*}
\begin{align*}
I_{2}&=\int_{0}^{t_{2}}(t_{2}-s)^{\alpha-1}E_{\alpha,\alpha}(t_{2}-s)\sigma(u(s))dW(s)\\
&\hspace{2mm}-\int_{0}^{t_{1}}(t_{1}-s)^{\alpha-1}E_{\alpha,\alpha}(t_{1}-s)\sigma(u(s))dW(s)\\
&=\int_{0}^{t_{1}}(t_{1}-s)^{\alpha-1}[E_{\alpha,\alpha}(t_{2}-s)-E_{\alpha,\alpha}(t_{1}-s)]\sigma (u(s))dW(s)\\
&\hspace{2mm}+\int_{0}^{t_{1}}[(t_{2}-s)^{\alpha-1}-(t_{1}-s)^{\alpha-1}]E_{\alpha,\alpha}(t_{2}-s)\sigma (u(s))dW(s)\\
&\hspace{2mm}+\int_{t_{1}}^{t_{2}}(t_{2}-s)^{\alpha-1}E_{\alpha,\alpha}(t_{2}-s)\sigma (u(s))dW(s)\\
&=: I_{21}+I_{22}+I_{23}.\tag{2.17}
\end{align*}

For any $0\leq \nu\leq 1$, by Lemma 2.6, we have
\begin{align*}
\mathbf{E}\|I_{1}\|_{\nu}^{2}=\mathbf{E}\|A^{\frac{\nu}{2}}[E_{\alpha}(t_{2})-E_{\alpha}(t_{1})]u_{0}\|^{2} \leq C_{\nu1}(t_{2}-t_{1})^{\nu\alpha}\mathbf{E}\|u_{0}\|^{2}.\tag{2.18}
\end{align*}

For the first term $I_{21}$ in (2.17), making use of Lemma 2.1 and Lemma 2.6 and Theorem 2.1, we get
\begin{align*}
\mathbf{E}\|I_{21}\|_{{\nu}}^{2}&=\mathbf{E}\|\int_{0}^{t_{1}}(t_{1}-s)^{\alpha-1}A^{\frac{\nu}{2}}[E_{\alpha,\alpha}(t_{2}-s)-E_{\alpha,\alpha}(t_{1}-s)]\sigma (u(s))dW(s)\|^{2}\\
&=\int_{0}^{t_{1}}\mathbf{E}\|(t_{1}-s)^{\alpha-1}A^{\frac{\nu}{2}}[E_{\alpha,\alpha}(t_{2}-s)-E_{\alpha,\alpha}(t_{1}-s)]\sigma (u(s))\|_{L_{2}^{0}}^{2}ds\\
&\leq C_{\nu2}^{2}C_{\sigma}^{2}(t_{2}-t_{1})^{\nu\alpha}(\int_{0}^{t_{1}}(t_{1}-s)^{2\alpha-2}\mathbf{E}\|u(s)\|^{2}ds)\\
&\leq \frac{C_{\nu2}^{2}C_{\sigma}^{2} T^{2\alpha-1}}{(2\alpha-1)}(t_{2}-t_{1})^{\nu\alpha}(\mathbf{E}\|u_{0}\|^{2}).\tag{2.19}
\end{align*}

Using Lemma 2.1, Lemma 2.4, Lemma 2.5 and Theorem 2.1, the term $I_{22}$ can be estimated as
\begin{align*}
\mathbf{E}\|I_{22}\|_{{\nu}}^{2}&=\mathbf{E}\|\int_{0}^{t_{1}}[(t_{2}-s)^{\alpha-1}-(t_{1}-s)^{\alpha-1}]A^{\frac{\nu}{2}}E_{\alpha,\alpha}(t_{2}-s)\sigma (u(s))dW(s)\|^{2}\\
&=\int_{0}^{t_{1}}\mathbf{E}\|[(t_{2}-s)^{\alpha-1}-(t_{1}-s)^{\alpha-1}]A^{\frac{\nu}{2}}E_{\alpha,\alpha}(t_{2}-s)\sigma (u(s))\|_{L_{2}^{0}}^{2}ds\\
&\leq C_{\alpha2}^{2}C_{\sigma}^{2}(\int_{0}^{t_{1}}[(t_{2}-s)^{\alpha-1}-(t_{1}-s)^{\alpha-1}]^{2}\mathbf{E}\|u(s)\|_{\nu}^{2}ds)\\
&\leq \frac{C_{\alpha2}^{2}C_{\sigma}^{2}T}{T_{0}^{2(1-\alpha)}}(t_{2}-t_{1})^{2(1-\alpha)}(\mathbf{E}\|u_{0}\|^{2}),\tag{2.20}
\end{align*}
where $0<T_{0}<T$.

For the term $I_{23}$, by Lemma 2.1, Lemma 2.4 and Theorem 2.1, there holds
\begin{align*}
\mathbf{E}\|I_{23}\|_{\nu}^{2}&=\mathbf{E}\|\int_{t_{1}}^{t_{2}}(t_{2}-s)^{\alpha-1}A^{\frac{\nu}{2}}E_{\alpha,\alpha}(t_{2}-s)\sigma (u(s))dW(s)\|^{2}\\
&=\int_{t_{1}}^{t_{2}}\mathbf{E}\|(t_{2}-s)^{\alpha-1}A^{\frac{\nu}{2}}E_{\alpha,\alpha}(t_{2}-s)\sigma (u(s))\|_{L_{2}^{0}}^{2}ds\\
&\leq C_{\alpha2}^{2}C_{\sigma}^{2} (\int_{t_{1}}^{t_{2}}(t_{2}-s)^{2\alpha-2}\mathbf{E}\|u(s)\|_{\nu}^{2}ds)\\
&\leq \frac{C_{\alpha2}^{2}C_{\sigma}^{2}}{(2\alpha-1)}(t_{2}-t_{1})^{2\alpha-1}(\mathbf{E}\|u_{0}\|^{2}).\tag{2.21}
\end{align*}

Taking expectation on both side of (2.16), and combining the estimates (2.18)-(2.21), we have
\begin{align*}
\mathbf{E}\|u(t_{2})-u(t_{1})\|_{\nu}^{2} \leq 2\mathbf{E}\|I_{1}\|_{\nu}^{2}+2\mathbf{E}\|I_{2}\|_{\nu}^{2}\leq C(t_{2}-t_{1})^{\beta},\tag{2.22}
\end{align*}
where $\beta=\max\{\nu\alpha,2(1-\alpha),(2\alpha-1)\}>0$.

This completes the proof of Theorem 2.2. $\square$

\section{Finite element method}

Let $\{\mathcal{T}_{h}\}_{h>0}$ be a family of regular triangulations of $D$ with the maximal mesh size of $h$. Denote by $V_{h}$ the spaces of continuous functions on $D$, which are piecewise polynomials with respect to $\mathcal{T}_{h}$. We define the discrete version of Laplacian to be an operator $A_{h}:V_{h}\rightarrow V_{h}$, it satisfies
\begin{align*}
(A_{h}\psi,\chi)=(\nabla\psi,\nabla\chi), \forall \psi,\chi \in V_{h}. \tag{3.1}
\end{align*}

The projection operator $P_{h}$ is the standard $L_{2}$-projection operator onto $V_{h}$, which defined by
\begin{align*}
(P_{h}v,\chi)=(v,\chi), \forall \chi \in V_{h}, \tag{3.2}
\end{align*}
where $v\in L_{2}(D)$.

It is easily shown that the operator $P_{h}v=\sum\limits_{j=0}^{N_{h}}a_{j}\varphi_{j}$ with a basis $\{\varphi_{i}\}_{i=0}^{N_{h}}$ can be solved from the equations $(\sum\limits_{j=0}^{N_{h}}a^{j}\varphi_{j},\varphi_{i})=(v,\varphi_{i})$.

By the definition of (3.1) we have
\begin{align*}
\|A_{h}^{\frac{1}{2}}\psi\|=\|\nabla\psi\|=\|A^{\frac{1}{2}}\psi\|=\|\psi\|_{1},\forall \psi \in V_{h}. \tag{3.3}
\end{align*}

Note that $P_{h}$ can be extended to $\dot{H}^{-1}$, that is, for all $\psi\in\dot{H}^{-1}$, there holds
\begin{align*}
\|A_{h}^{-\frac{1}{2}}P_{h}\psi\|&=\sup\limits_{\varphi\in V_{h}}\frac{(\psi,\varphi)}{\|A_{h}^{\frac{1}{2}}\varphi\|}=\sup\limits_{\varphi\in V_{h}}\frac{(\psi,\varphi)}{\|A^{\frac{1}{2}}\varphi\|}\leq \sup\limits_{\varphi\in \dot{H}^{-1}}\frac{(\psi,\varphi)}{\|A^{\frac{1}{2}}\varphi\|}\\
&=\|A^{-\frac{1}{2}}\psi\|.\tag{3.4}
\end{align*}

From (3.3) and (3.4), it is easy to get
\begin{align*}
\|A_{h}^{\frac{1}{2}}P_{h}\psi\|\leq C\|A^{\frac{1}{2}}\psi\|=C\|\psi\|_{1}, \forall \psi\in\dot{H}^{1}. \tag{3.5}
\end{align*}

Interpolation between (3.4) and (3.5) yields
\begin{align*}
\|A_{h}^{\rho}P_{h}\psi\|\leq C\|A^{\rho}\psi\|, \forall \psi\in\dot{H}^{\rho},\rho\in[-\frac{1}{2},\frac{1}{2}].  \tag{3.6}
\end{align*}

Next, we will present a standard Galerkin finite element method for time-fractional stochastic heat equation (2.2).

\subsection{Semidiscrete finite element approximation}

The semi-discretized version of (2.2) is to find a process $u_{h}(t)=u_{h}(.,t)\in V_{h}$ such that
\begin{align*}
\begin{cases}
^{C}D_{t}^{\alpha}u_{h}(t)=A_{h}u_{h}(t)+P_{h}\sigma(u_{h}(t))\frac{d W(t)}{d t},\\
u_{h}(0)=P_{h}u_{0}=u_{h0}. \end{cases}\tag{3.7}
\end{align*}

The mild solution of Galerkin approximation $u_{h}(t):[0,T]\times\Omega\rightarrow V_{h}$ in (3.7) is given by
\begin{align*}
u_{h}(t)=E_{\alpha}^{h}(t)P_{h}u_{0}+\int_{0}^{t}(t-s)^{\alpha-1}E_{\alpha,\alpha}^{h}(t-s)P_{h}\sigma(u_{h}(s))dW(s),\tag{3.8}
\end{align*}
in which
\begin{align*}
E_{\alpha}^{h}(t)=\int_{0}^{\infty}\xi_{\alpha}(\theta)S_{h}(t^{\alpha}\theta)d\theta,~E_{\alpha,\alpha}^{h}(t)=\int_{0}^{\infty}\alpha\theta\xi_{\alpha}(\theta)S_{h}(t^{\alpha}\theta)d\theta,
\end{align*}
where $S_{h}(t)=e^{t A_{h}}$ is the analytic semigroup generated by $A_{h}$.

To prove our main results, some useful results of the corresponding deterministic problem will be needed below.

\textbf{Lemma 3.1.} (see [14]) Let $\rho\geq0$, there exists a constant $C_{\rho}>0$ such that
\begin{align*}
\|A_{h}^{\rho}S_{h}(t)v\|\leq C_{\rho}t^{-\rho}\|v\|, ~\forall t>0,h\in(0,1). \tag{3.9}
\end{align*}

\textbf{Lemma 3.2.} (see [14]) Let $F_{h}(t)=S_{h}(t)P_{h}-S(t),t\geq0$, for $0\leq\nu\leq\mu\leq2$, there exists a constant $C>0$ such that
\begin{align*}
\|F_{h}(t)v\|\leq Ch^{\mu}t^{-\frac{\mu-\nu}{2}}\|v\|_{\nu}, ~\forall v\in \dot{H}^{\nu},t>0,h\in(0,1).\tag{3.10}
\end{align*}

\textbf{Lemma 3.3.} For $0\leq\mu<1$, there exist two constants $C_{\mu1}>0$ and $C_{\mu2}>0$ such that
\begin{align*}
\|A_{h}^{\mu}E_{\alpha}^{h}(t)P_{h}v\|\leq C_{\mu1}t^{-\mu\alpha}\|v\|,~\|A_{h}^{\mu}E_{\alpha,\alpha}^{h}(t)P_{h}v\|\leq C_{\mu2}t^{-\mu\alpha}\|v\|.\tag{3.11}
\end{align*}

\textbf{Proof.} For any $0\leq\mu<1$, by using Lemma 2.2, Lemma 3.1 and the inequality (3.6) ($\rho=0$), we have
\begin{align*}
\|A_{h}^{\mu}E_{\alpha}^{h}(t)P_{h}v\|&=\|\int_{0}^{\infty}\xi_{\alpha}(\theta)A_{h}^{\mu}S_{h}(t^{\alpha}\theta)P_{h}vd\theta\|\\
&\leq C_{\mu}t^{-\mu\alpha}(\int_{0}^{\infty}\theta^{-\mu}\xi_{\alpha}(\theta)\|v\|d\theta)\\
&=\frac{C_{\mu}\Gamma(1-\mu)}{\Gamma(1-\alpha\mu)}t^{-\mu\alpha}\|v\|,
\end{align*}
and
\begin{align*}
\|A_{h}^{\mu}E_{\alpha,\alpha}^{h}(t)P_{h}v\|&=\|\int_{0}^{\infty}\alpha\theta\xi_{\alpha}(\theta)A_{h}^{\mu}S_{h}(t^{\alpha}\theta)P_{h}vd\theta\|\\
&\leq C_{\mu}\alpha t^{-\mu\alpha}(\int_{0}^{\infty}\theta^{1-\mu}\xi_{\alpha}(\theta)\|v\|d\theta)\\
&=\frac{C_{\mu}\alpha \Gamma(2-\mu)}{\Gamma(1+\alpha(1-\mu))}t^{-\mu\alpha}\|v\|.
\end{align*}

\textbf{Lemma 3.4.} For any $0\leq r<1$, there exist two constants $C_{r1}>0$ and $C_{r2}>0$ such that
\begin{align*}
\|[E_{\alpha}^{h}(t)P_{h}-E_{\alpha}(t)]v\|\leq C_{r1}t^{-\frac{r\alpha}{2}}h^{1+r}\|v\|_{1},\tag{3.12}
\end{align*}
and
\begin{align*}
\|[E_{\alpha,\alpha}^{h}(t)P_{h}-E_{\alpha,\alpha}(t)]v\|\leq C_{r2}t^{-\frac{r\alpha}{2}}h^{1+r}\|v\|_{1}.\tag{3.13}
\end{align*}

\textbf{Proof.} Using Lemma 2.2 and Lemma 3.2, setting $\mu=1+r,r\in[0,1)$ and $\nu=1$ in (3.10), we deduce that
\begin{align*}
\|[E_{\alpha}^{h}(t)P_{h}-E_{\alpha}(t)]v\|&=\|\int_{0}^{\infty}\xi_{\alpha}(\theta)F_{h}(t^{\alpha}\theta)vd\theta\|\\
&\leq \int_{0}^{\infty}\xi_{\alpha}(\theta)\|F_{h}(t^{\alpha}\theta)v\|d\theta\\
&\leq Ct^{-\frac{r\alpha}{2}} h^{1+r}(\int_{0}^{\infty}\theta^{-\frac{r}{2}}\xi_{\alpha}(\theta)\|v\|_{1}d\theta)\\
&=\frac{C\Gamma(1-\frac{r}{2})}{\Gamma(1-\frac{r\alpha}{2})}t^{-\frac{r\alpha}{2}}h^{1+r}\|v\|_{1},
\end{align*}
and
\begin{align*}
\|[E_{\alpha,\alpha}^{h}(t)P_{h}-E_{\alpha,\alpha}(t)]v\|&=\|\int_{0}^{\infty}\alpha\theta\xi_{\alpha}(\theta)F_{h}(t^{\alpha}\theta)vd\theta\|\\
&\leq\int_{0}^{\infty}\alpha\theta\xi_{\alpha}(\theta)\|F_{h}(t^{\alpha}\theta)v\|d\theta\\
&\leq C\alpha t^{-\frac{r\alpha}{2}} h^{1+r}(\int_{0}^{\infty}\theta^{1-\frac{r}{2}}\xi_{\alpha}(\theta)\|v\|_{1}d\theta)\\
&=\frac{C\alpha\Gamma(2-\frac{r}{2})}{\Gamma(1+\frac{\alpha(2-r)}{2})}t^{-\frac{r\alpha}{2}}h^{1+r}\|v\|_{1}.
\end{align*}

\textbf{Remark 3.1.} If we take $\mu=1+r,r\in[0,1)$ and $\nu=0$ in (3.10), then the following inequalities hold:
\begin{align*}
\|[E_{\alpha}^{h}(t)P_{h}-E_{\alpha}(t)]v\|\leq C_{r1}t^{-\frac{(1+r)\alpha}{2}}h^{1+r}\|v\|,
\end{align*}
and
\begin{align*}
\|[E_{\alpha,\alpha}^{h}(t)P_{h}-E_{\alpha,\alpha}(t)]v\|\leq C_{r2}t^{-\frac{(1+r)\alpha}{2}}h^{1+r}\|v\|.
\end{align*}

\textbf{Lemma 3.5.} Assume that $\sigma$ satisfies (2.3), for any $t\in [0,T]$ and $\alpha\in(\frac{1}{2},1)$. Let $u_{h}(t)$ be a mild solution to (3.7). There exists a constant $C>0$ such that
\begin{align*}
\sup\limits_{t\in [0,T]}\|u_{h}(t)\|_{L_{2}(\Omega,H)}^{2}\leq C\|u_{0}\|_{L_{2}(\Omega,H)}^{2}. \tag{3.14}
\end{align*}
\textbf{Proof.} For any $t\in [0,T]$, from (3.8) and based on Lemma 2.1 and Lemma 3.3 ($\mu$=0), we obtain
\begin{align*}
\mathbf{E}\|u_{h}(t)\|^{2}&\leq 2\mathbf{E}\|E_{\alpha}^{h}(t)P_{h}u_{0}\|^{2}+2\mathbf{E}\|\int_{0}^{t}(t-s)^{\alpha-1}E_{\alpha,\alpha}^{h}(t-s)P_{h}\sigma(u_{h}(s))dW(s)\|^{2}\\
&\leq 2C_{\mu1}^{2}\mathbf{E}\|u_{0}\|^{2}+2(\int_{0}^{t}\mathbf{E}\|(t-s)^{\alpha-1}E_{\alpha,\alpha}^{h}(t-s)P_{h}\sigma(u_{h}(s))\|_{L_{2}^{0}}^{2}ds)\\
&\leq 2C_{\mu1}^{2} \mathbf{E}\|u_{0}\|^{2}+2C_{\mu2}^{2}C_{\sigma}^{2}(\int_{0}^{t}(t-s)^{(2\alpha-1)-1}\mathbf{E}\|u_{h}(s)\|^{2}ds).
\end{align*}

Thus, by virtue of Lemma 2.7, we deduce that
\begin{align*}
\sup\limits_{t\in [0,T]}\mathbf{E}\|u_{h}(t)\|^{2}\leq C\mathbf{E}\|u_{0}\|^{2}.
\end{align*}

The proof of the lemma is completed.  $\square$

\textbf{Theorem 3.1.} For all $t\in[0,T]$, $\alpha\in(\frac{1}{2},1)$ and $r\in[0,1)$, let $u_{h}(t)$ and $u(t)$ be the mild solutions to (3.7) and (2.2), respectively. Then there exist a constant $C>0$, which
is independent of $h$, such that
\begin{align*}
\|u_{h}(t)-u(t)\|_{L_{2}(\Omega,H)}^{2}\leq Ch^{2(1+r)}.
\end{align*}

\textbf{Proof.} For any $t\in[0,T]$, from (3.8) and (2.4), we have
\begin{align*}
u_{h}(t)-u(t)&=E_{\alpha}^{h}(t)P_{h}u_{0}-E_{\alpha}(t)u_{0}\\
&\hspace{2mm}+\int_{0}^{t}(t-s)^{\alpha-1}E_{\alpha,\alpha}^{h}(t-s)P_{h}\sigma(u_{h}(s))dW(s)\\
&\hspace{2mm}-\int_{0}^{t}(t-s)^{\alpha-1}E_{\alpha,\alpha}(t-s)\sigma(u(s))dW(s)\\
&:=J_{1}+J_{2},\tag{3.15}
\end{align*}
where
\begin{align*}
J_{1}=[E_{\alpha}^{h}(t)P_{h}-E_{\alpha}(t)]u_{0},
\end{align*}
\begin{align*}
J_{2}&=\int_{0}^{t}(t-s)^{\alpha-1}[E_{\alpha,\alpha}^{h}(t-s)P_{h}\sigma(u_{h}(s))-E_{\alpha,\alpha}(t-s)\sigma(u(s))]dW(s)\\
&=\int_{0}^{t}(t-s)^{\alpha-1}E_{\alpha,\alpha}^{h}(t-s)P_{h}[\sigma(u_{h}(s))-\sigma(u(s))]dW(s)\\
&\hspace{2mm}+\int_{0}^{t}(t-s)^{\alpha-1}[E_{\alpha,\alpha}^{h}(t-s)P_{h}-E_{\alpha,\alpha}(t-s)](\sigma(u(s))-\sigma(u(t)))dW(s)\\
&\hspace{2mm}+\int_{0}^{t}(t-s)^{\alpha-1}[E_{\alpha,\alpha}^{h}(t-s)P_{h}-E_{\alpha,\alpha}(t-s)]\sigma(u(t))dW(s)\\
&:=J_{21}+J_{22}+J_{23}.\tag{3.16}
\end{align*}

The application of Lemma 3.4, the term $J_{1}$ can be estimated as
\begin{align*}
\mathbf{E}\|J_{1}\|=\mathbf{E}\|[E_{\alpha}^{h}(t)P_{h}-E_{\alpha}(t)]u_{0}\|\leq C_{r1}t^{-\frac{r\alpha}{2}}h^{1+r}\mathbf{E}\|u_{0}\|_{1}. \tag{3.17}
\end{align*}

For the term $J_{21}$ in (3.16), applying Lemma 2.1, Lemma 3.3 ($\mu=0$) and the Lipschitz condition (2.3), we obtain
\begin{align*}
\mathbf{E}\|J_{21}\|^{2}&=\mathbf{E}\|\int_{0}^{t}(t-s)^{\alpha-1}E_{\alpha,\alpha}^{h}(t-s)P_{h}[\sigma(u_{h}(s))-\sigma(u(s))]dW(s)\|^{2}\\
&=\int_{0}^{t}\mathbf{E}\|(t-s)^{\alpha-1}E_{\alpha,\alpha}^{h}(t-s)P_{h}[\sigma(u_{h}(s))-\sigma(u(s))]\|_{L_{2}^{0}}^{2}ds\\
&\leq C\int_{0}^{t}(t-s)^{2\alpha-2}\mathbf{E}\|u_{h}(s)-u(s)\|^{2}ds.  \tag{3.18}
\end{align*}

A combination of Lemma 2.1, Lemma 3.4, and Theorem 2.2, we deduce that
\begin{align*}
\mathbf{E}\|J_{22}\|^{2}&=\mathbf{E}\|\int_{0}^{t}(t-s)^{\alpha-1}[E_{\alpha,\alpha}^{h}(t-s)P_{h}-E_{\alpha,\alpha}(t-s)](\sigma(u(s))-\sigma(u(t)))dW(s)\|^{2}\\
&=\int_{0}^{t}\mathbf{E}\|(t-s)^{\alpha-1}[E_{\alpha,\alpha}^{h}(t-s)P_{h}-E_{\alpha,\alpha}(t-s)](\sigma(u(s))-\sigma(u(t)))\|_{L_{2}^{0}}^{2}ds\\
&\leq Ch^{2(1+r)}(\int_{0}^{t}(t-s)^{2\alpha-2-r\alpha}\mathbf{E}\|u(t)-u(s)\|_{1}^{2}ds)\\
&\leq Ch^{2(1+r)}(\int_{0}^{t}(t-s)^{2\alpha-2-r\alpha+\beta}ds)\\
&\leq\frac{CT^{(2-r)\alpha+\beta-1}}{((2-r)\alpha+\beta-1)}h^{2(1+r)}, \tag{3.19}
\end{align*}
where the  parameters $r$ and $\alpha$ should satisfy $(2-r)\alpha+\beta-1>0$.

Similarly, for the term $J_{23}$, by Theorem 2.1, we get
\begin{align*}
\mathbf{E}\|J_{23}\|^{2}&=\mathbf{E}\|\int_{0}^{t}(t-s)^{\alpha-1}[E_{\alpha,\alpha}^{h}(t-s)P_{h}-E_{\alpha,\alpha}(t-s)]\sigma(u(t))dW(s)\|^{2}\\
&=\int_{0}^{t}\mathbf{E}\|(t-s)^{\alpha-1}[E_{\alpha,\alpha}^{h}(t-s)P_{h}-E_{\alpha,\alpha}(t-s)]\sigma(u(t))\|_{L_{2}^{0}}^{2}ds\\
&\leq Ch^{2(1+r)}(\int_{0}^{t}(t-s)^{2\alpha-2-r\alpha}\mathbf{E}\|u(t)\|_{1}^{2}ds)\\
&\leq \frac{CT^{(2-r)\alpha-1}}{((2-r)\alpha-1)}(\mathbf{E}\|u_{0}\|^{2})h^{2(1+r)},\tag{3.20}
\end{align*}
where the parameter $r$ should satisfy $(2-r)\alpha-1>0$.

Taking expectation on (3.15) and together with (3.17)-(3.20), we have
\begin{align*}
\mathbf{E}\|u_{h}(t)-u(t)\|^{2}&\leq 2\mathbf{E}\|J_{1}\|^{2} +2\mathbf{E}\|J_{2}\|^{2} \\
&\leq C_{1}h^{2(1+r)}+C_{2}\int_{0}^{t}(t-s)^{(2\alpha-1)-1}\mathbf{E}\|u_{h}(s)-u(s)\|^{2}ds.\tag{3.21}
\end{align*}

Therefore, by means of Lemma 2.7, there holds
\begin{align*}
\mathbf{E}\|u_{h}(t)-u(t)\|^{2}\leq C(C_{1},C_{2},T,\alpha)h^{2(1+r)}.
\end{align*}

The proof is completed.   $\square$ 

\subsection{Fully discrete schemes}

Denote by $t_{n}=n\tau,n=0,1,\cdots, N$ the time mesh point with a fixed time mesh size $\tau>0$, which satisfy the integration time $0\leq t_{n}\leq T$ and $\tau=\frac{T}{N}$.  Then the semi-discretized version of mild solution (3.8) at time $t_{n}$ is shown that
\begin{align*}
u_{h}(t_{n})=E_{\alpha}^{h}(t_{n})P_{h}u_{0}+\int_{0}^{t_{n}}(t_{n}-s)^{\alpha-1}E_{\alpha,\alpha}^{h}(t_{n}-s)P_{h}\sigma(u_{h}(s))dW(s).\tag{3.22}
\end{align*}

Now we will introduce the Gorenflo-Mainardi-Moretti-Paradisi (GMMP) scheme, which was firstly developed in [15]. Then the Caputo fractional derivative can be approximated by
\begin{align*}
^{C}D_{t}^{\alpha}u(x,t_{n})&\approx \frac{1}{\tau^{\alpha}}\sum\limits_{k=0}^{n}\omega_{k}^{\alpha}[u(x,t_{n-k})-u(x,0)]\\
&=\frac{1}{\tau^{\alpha}}[\sum\limits_{k=0}^{n}\omega_{k}^{\alpha}u(x,t_{n-k})-b_{n}u(x,0)],  \tag{3.23}
\end{align*}
where
\begin{align*}
\omega_{k}^{\alpha}=(-1)^{k}\dbinom{\alpha}{k}=\frac{\Gamma(k-\alpha)}{\Gamma(-\alpha)\Gamma(k+1)},
\end{align*}
and
\begin{align*}
b_{n}=\sum\limits_{k=0}^{n}\omega_{k}^{\alpha}=\frac{\Gamma(n+1-\alpha)}{\Gamma(1-\alpha)\Gamma(n+1)},n\geq0.
\end{align*}

Furthermore, $\omega_{k}^{\alpha}$ and $b_{n}$ have the following properties.

\textbf{Lemma 3.6.} (see [16,17]) For $\alpha>0$, $n=1,2,\cdots$, we have

(1) $\omega_{0}^{\alpha}=1$, $\omega_{n}^{\alpha}<0$, $|\omega_{n+1}^{\alpha}|<|\omega_{n}^{\alpha}|$, and  $0<-\sum\limits_{k=1}^{n}\omega_{k}^{\alpha}<-\sum\limits_{k=1}^{\infty}\omega_{k}^{\alpha}=\omega_{0}^{\alpha}$.

(2) $b_{n}-b_{n-1}=\omega_{n}^{\alpha}<0$, i.e., $b_{n}< b_{n-1}< b_{n-2}<\cdots< b_{0}=1$.

By using the GMMP scheme (3.23), we denote $u^{n}\approx u(t_{n})$ as the approximation of $u(t_{n})$. Then the full discrete scheme for equation (2.2) can be defined by seeking an $\mathcal{F}_{t_{n}}$-adapted process ${u_{h}^{n}}$ satisfying:
\begin{align*}
\begin{cases}
\frac{1}{\tau^{\alpha}}[\sum\limits_{k=0}^{n}\omega_{k}^{\alpha}u_{h}^{n-k}-b_{n}u_{h}^{0}]=A_{h}u_{h}^{n}+\frac{1}{\tau}\int_{t_{n-1}}^{t_{n}}P_{h}\sigma(u_{h}^{n-1})d W(s), \\
u_{h}^{0}=P_{h}u_{0}.
\end{cases}     \tag{3.24}
\end{align*}

With the definition of $R(\lambda,X)=(\lambda I-X)^{-1},\lambda>0$, and $E_{\tau h}=R(\tau^{-\alpha},A_{h})=(\tau^{-\alpha}I-A_{h})^{-1}$. The above scheme (3.24) can be rewritten as:
\begin{align*}
\begin{cases}
u_{h}^{n}=\tau^{-\alpha}b_{n}E_{\tau h}u_{h}^{0}-\tau^{-\alpha}E_{\tau h}\sum\limits_{k=1}^{n}\omega_{k}^{\alpha}u_{h}^{n-k}+\frac{1}{\tau}\int_{t_{n-1}}^{t_{n}}E_{\tau h}P_{h}\sigma(u_{h}^{n-1})d W(s),\\
u_{h}^{0}=P_{h}u_{0}.\end{cases} \tag{3.25}
\end{align*}

\textbf{Lemma 3.7.} For any $\tau>0$ and $h\in(0,1)$. There exists a constant $C>0$ such that
\begin{align*}
\|E_{\tau h}v\|\leq C\tau^{\alpha}\|v\|,~\|E_{\tau h}P_{h}v\|\leq C\tau^{\alpha}\|v\|, ~ \forall~ v\in H.\tag{3.26}
\end{align*}

\textbf{Proof.} For $\lambda>0$, applying Lemma 3.1 ($\rho=0$), we get
\begin{align*}
\|R(\lambda,A_{h})v\|&=\|(\lambda I-A_{h})^{-1}v\|=\|\int_{0}^{\infty}e^{-(\lambda I-A_{h}) t}vdt\|\\
&=\|\int_{0}^{\infty}e^{-\lambda t}S_{h}(t)vdt\|\leq\int_{0}^{\infty}e^{-\lambda t}\|S_{h}(t)v\|dt\\
&\leq \frac{C}{\lambda}\|v\|.
\end{align*}

Therefore, we have $\|E_{\tau h}v\|=\|R(\tau^{-\alpha},A_{h})v\|\leq C\tau^{\alpha}\|v\|$.

Using the inequality (3.6) ($\rho=0$) and Lemma 3.1 again, we obtain
\begin{align*}
\|R(\lambda,A_{h})P_{h}v\|&=\|(\lambda I-A_{h})^{-1}P_{h}v\|= \|\int_{0}^{\infty}e^{-(\lambda I-A_{h}) t}P_{h}vdt\|\\
&= \|\int_{0}^{\infty}e^{-\lambda t}S_{h}(t)P_{h}vdt\|\\
&\leq\int_{0}^{\infty}e^{-\lambda t}\|S_{h}(t)P_{h}v\|dt\\
&\leq \frac{C}{\lambda}\|v\|.
\end{align*}

Based on above estimate, it is easy to get $\|E_{\tau h}P_{h}v\|=\|R(\tau^{-\alpha},A_{h})P_{h}v\|\leq C\tau^{\alpha}\|v\|$.

\textbf{Lemma 3.8.} For any $\lambda>0$ and $\mu\in R$, there exists a constant $C$ such that
\begin{align*}
\|[\mu R(\lambda,A_{h})-I]P_{h}v\|\leq C\lambda^{-1}\|v\|. \tag{3.27}
\end{align*}

\textbf{Proof.} For any $\lambda>0$ and $\mu\in R$, using the same argument as the proof of (3.26), based on the estimate $\|A_{h}^{\rho}S_{h}(t)v\|\leq C_{\rho}T_{0}^{-\rho}\|v\|, 0<T_{0}\leq t<\infty$, one can deduce that
\begin{align*}
&\|[\mu R(\lambda,A_{h})-I]P_{h}v\|\\
&=\|[\mu(\lambda I-A_{h})^{-1}-(\lambda I-A_{h})^{-1}(\lambda I-A_{h})]P_{h}v\|\\
&=\|(\lambda I-A_{h})^{-1}[(\mu-\lambda)I+A_{h}]P_{h}v\|\\
&\leq \|(|\mu|+|\lambda|)(\lambda I-A_{h})^{-1}P_{h}v\|+\|(\lambda I-A_{h})^{-1}A_{h}P_{h}v\|\\
&\leq\frac{C}{\lambda}\|v\|+\int_{0}^{\infty}e^{-\lambda t}\|A_{h}S_{h}(t)P_{h}v\|dt\\
&\leq\frac{C}{\lambda}\|v\|+\frac{CT_{0}^{-1}}{\lambda}\|v\|\\
&\leq C^{\dagger}\lambda^{-1}\|v\|.
\end{align*}

\textbf{Lemma 3.9.} For any $t\in[0,T]$, there exists a constant $C>0$ such that
\begin{align*}
\|[I-S_{h}(t)]P_{h}v\|\leq C t\|A_{h}v\|. \tag{3.28}
\end{align*}

\textbf{Proof.} Using Lemma 3.1 and the estimate (3.6) we have
\begin{align*}
\|[I-S_{h}(t)]P_{h}v\|&=\|(-\int_{0}^{t}A_{h}e^{sA_{h}}ds)P_{h}v\|\\
&\leq\int_{0}^{t}\|S_{h}(s)A_{h}P_{h}v\|ds\\
&\leq C(\int_{0}^{t}\|A_{h}v\|dt)\\
&=Ct\|A_{h}v\|.
\end{align*}

Now we shall provide the error estimates of the numerical solution $u_{h}^{n}$ to $u(t_{n})$. Denote by $e^{n}=u_{h}^{n}-u(t_{n})$, then we have the following results.

\textbf{Theorem 3.2.} For any $\alpha\in(\frac{1}{2},1)$ and $r\in [0,1)$, let $u_{h}^{n}$ and $u(t_{n})$ be the solutions to (3.25) and (2.2), respectively. Then there exist a constant $C>0$ such that
\begin{align*}
\|e^{n}\|_{L_{2}(\Omega,H)}^{2}\leq C[\tau^{2\alpha}+h^{2(1+r)}].\tag{3.29}
\end{align*}

\textbf{Proof.} We rewrite $e^{n}=u_{h}^{n}-u(t_{n})=[u_{h}^{n}-u_{h}(t_{n})]+[u_{h}(t_{n})-u(t_{n})]:=\xi^{n}+\eta^{n}$. The estimate for $\eta^{n}$ can be obtained by Theorem 3.1, that is,
\begin{align*}
\|\eta^{n}\|_{L_{2}(\Omega,H)}^{2}\leq Ch^{2(1+r)}.\tag{3.30}
\end{align*}

To estimate the term $\xi^{n}$, making use of (3.25) and (3.22), we obtain
\begin{align*}
\xi^{n}&=[\tau^{-\alpha}b_{n}E_{\tau h}u_{h}^{0}-E_{\alpha}^{h}(t_{n})P_{h}u_{0}]-\tau^{-\alpha}E_{\tau h}\sum\limits_{k=1}^{n}\omega_{k}^{\alpha}u_{h}^{n-k}\\
&\hspace{2mm}+\frac{1}{\tau}\int_{t_{n-1}}^{t_{n}}E_{\tau h}P_{h}\sigma(u_{h}^{n-1})d W(s)\\
&\hspace{2mm}-\int_{0}^{t_{n}}(t_{n}-s)^{\alpha-1}E_{\alpha,\alpha}^{h}(t_{n}-s)P_{h}\sigma(u_{h}(s))dW(s)\\
&:=L_{1}+L_{2}+L_{3}+L_{4}. \tag{3.31}
\end{align*}

The term $L_{1}$ can be rewritten as
\begin{align*}
L_{1}&=\tau^{-\alpha}b_{n}E_{\tau h}u_{h}^{0}-E_{\alpha}^{h}(t_{n})P_{h}u_{0}\\
&=[\tau^{-\alpha}b_{n}E_{\tau h}u_{h}^{0}-P_{h}u_{0}]+[P_{h}u_{0}-E_{\alpha}^{h}(t_{n})P_{h}u_{0}]\\
&:=L_{11}+L_{12}.  \tag{3.32}
\end{align*}

Setting $\mu=\tau^{-\alpha}b_{n}$, by using Lemma 3.8, we have
\begin{align*}
\mathbf{E}\|L_{11}\|^{2}=\mathbf{E}\|\mu E_{\tau h}u_{h}^{0}-P_{h}u_{0}\|^{2}=\mathbf{E}\|[\mu R(\tau^{-\alpha},A_{h})-I]P_{h}u_{0}\|^{2}\leq C\tau^{2\alpha}\mathbf{E}\|u_{0}\|^{2}.  \tag{3.33}
\end{align*}

Using Lemma 2.2 and Lemma 3.9, we obtain
\begin{align*}
\mathbf{E}\|L_{12}\|^{2}&=\mathbf{E}\|P_{h}u_{0}-E_{\alpha}^{h}(t_{n})P_{h}u_{0}\|^{2}\\
&=\mathbf{E}\|\int_{0}^{\infty}\xi_{\alpha}(\theta)P_{h}u_{0}d\theta-\int_{0}^{\infty}\xi_{\alpha}(\theta)S_{h}(t_{n}^{\alpha}\theta)P_{h}u_{0}d\theta\|^{2}\\
&=\mathbf{E}\|\int_{0}^{\infty}\xi_{\alpha}(\theta)[I-S_{h}(t_{n}^{\alpha}\theta)]P_{h}u_{0}d\theta\|^{2}\\
&\leq\int_{0}^{\infty}\xi_{\alpha}(\theta)\mathbf{E}\|[I-S_{h}(t_{n}^{\alpha}\theta)]P_{h}u_{0}\|^{2}d\theta\\
&\leq Ct_{n}^{2\alpha}(\int_{0}^{\infty}\theta\xi_{\alpha}(\theta)\mathbf{E}\|A_{h}u_{0}\|^{2}d\theta)\\
&= \frac{Cn^{2\alpha}\Gamma(2)}{\Gamma(1+\alpha)}\tau^{2\alpha}\mathbf{E}\|A_{h}u_{0}\|^{2}. \tag{3.34}
\end{align*}

According to Lemma 3.6, we have $|\omega_{k}^{\alpha}|^{2}<|\omega_{0}^{\alpha}|^{2}=1$, the application of Lemma 3.5 and Lemma 3.7 yields
\begin{align*}
\mathbf{E}\|L_{2}\|^{2}&=\mathbf{E}\|-\tau^{-\alpha}E_{\tau h}\sum\limits_{k=1}^{n}\omega_{k}^{\alpha}u_{h}^{n-k}\|^{2}\\
&=\mathbf{E}\|-\tau^{-\alpha}E_{\tau h}\sum\limits_{k=1}^{n}\omega_{k}^{\alpha}[(u_{h}^{n-k}-u_{h}(t_{n-k}))+u_{h}(t_{n-k})]\|^{2}\\
&\leq C\sum\limits_{k=1}^{n}\mathbf{E}\|\xi^{n-k}\|^{2}+C^{\dag}\sum\limits_{k=1}^{n}\mathbf{E}\|E_{\tau h}u_{h}(t_{n-k})\|^{2}\\
&\leq C\sum\limits_{k=1}^{n}\mathbf{E}\|\xi^{n-k}\|^{2}+ C^{\ddag}\tau^{2\alpha}(\mathbf{E}\|u_{0}\|^{2}). \tag{3.35}
\end{align*}

Making use of Lemma 2.1, Lemma 3.5, Lemma 3.7 and the condition (2.3), the estimate for the term $L_{3}$ in (3.31) gives
\begin{align*}
\mathbf{E}\|L_{3}\|^{2}&=\mathbf{E}\|\frac{1}{\tau}\int_{t_{n-1}}^{t_{n}}E_{\tau h}P_{h}\sigma(u_{h}^{n-1})d W(s)\|^{2}\\
&\leq 2\mathbf{E}\|\frac{1}{\tau}\int_{t_{n-1}}^{t_{n}}E_{\tau h}P_{h}[\sigma(u_{h}^{n-1})-\sigma(u_{h}(t_{n-1}))]dW(s)\|^{2}\\
&\hspace{2mm}+2\mathbf{E}\|\frac{1}{\tau}\int_{t_{n-1}}^{t_{n}}E_{\tau h}P_{h}\sigma(u_{h}(t_{n-1}))dW(s)\|^{2}\\
&= \frac{2}{\tau}\int_{t_{n-1}}^{t_{n}}\mathbf{E}\|E_{\tau h}P_{h}[\sigma(u_{h}^{n-1})-\sigma(u_{h}(t_{n-1}))]\|_{L_{2}^{0}}^{2}ds\\
&\hspace{2mm}+\frac{2}{\tau}\int_{t_{n-1}}^{t_{n}}\mathbf{E}\|E_{\tau h}P_{h}\sigma(u_{h}(t_{n-1}))\|_{L_{2}^{0}}^{2}ds\\
&\leq \frac{2C_{\sigma}^{2}}{\tau}\int_{t_{n-1}}^{t_{n}}\tau^{2\alpha}\mathbf{E}\|\xi^{n-1}\|^{2}ds+\frac{2C_{\sigma}^{2}}{\tau}\int_{t_{n-1}}^{t_{n}}\tau^{2\alpha}\mathbf{E}\|u_{h}(t_{n-1})\|^{2}ds\\
&\leq C^{\dag}\tau^{2\alpha}\mathbf{E}\|\xi^{n-1}\|^{2}+ C^{\ddag}\tau^{2\alpha}(\mathbf{E}\|u_{0}\|^{2}).\tag{3.36}
\end{align*}

By virtue of Lemma 2.1, Lemma 3.3, Lemma 3.5 and the condition (2.3), we have
\begin{align*}
\mathbf{E}\|L_{4}\|^{2}&=\mathbf{E}\|-\int_{0}^{t_{n}}(t_{n}-s)^{\alpha-1}E_{\alpha,\alpha}^{h}(t_{n}-s)P_{h}\sigma(u_{h}(s))dW(s)\|^{2}\\
&=\int_{0}^{t_{n}}\mathbf{E}\|(t_{n}-s)^{\alpha-1}E_{\alpha,\alpha}^{h}(t_{n}-s)P_{h}\sigma(u_{h}(s))\|_{L_{2}^{0}}^{2}ds\\
&\leq C(\int_{0}^{t_{n}}(t_{n}-s)^{2\alpha-2}\mathbf{E}\|u_{h}(s)\|^{2}ds)\\
&\leq \frac{Cn^{2\alpha}}{T_{0}(2\alpha-1)}\tau^{2\alpha}(\mathbf{E}\|u_{0}\|^{2}),\tag{3.37}
\end{align*}
where $0<T_{0}<T$.

Therefore, taking expectation on both side of (3.31), collecting all the above terms and applying a discrete version of Gronwall's Lemma, we have
\begin{align*}
\|\xi^{n}\|_{L_{2}(\Omega,H)}^{2}\leq C^{\dag}\tau^{2\alpha}.\tag{3.38}
\end{align*}

Thus, using the triangle inequality and combining (3.30) and (3.38), it gives (3.29). This completes the proof of Theorem 3.2.   $\square$

\section{Numerical example}

In this section, we focus on testing the achieved theoretical convergence results obtained in the previous section. For the numerical illustration of error bounds in Theorem 3.1 and 3.2. We consider the following one-dimension time-fractional stochastic heat equation:
\begin{align*}
\begin{cases}
 ^{C}D_{t}^{\alpha}u(x,t)=\frac{\partial^{2}u(x,t)}{\partial x^{2}}+u(x,t)\frac{dW(t)}{dt},~x\in (0,2),t\in(0,1],\\
u(x,0)=\sin x,~x\in (0,2),\\
u(0,t)=u(2,t)=0,~t\in (0,1].\\
\end{cases}
\end{align*}

To approximate the stochastic integral, we define a partition of $[0,T]$ by intervals $[t_{n-1},t_{n}]$ for $n=1,2,\ldots,N$, where $t_{n}=n\tau,\tau=T/N$. A sequence of noise which approximates the white noise is defined as
\begin{align*}
\frac{dW_{n}(t)}{dt}=\sum\limits_{n=1}^{N}\frac{1}{\sqrt{\tau}}\zeta_{n}e_{n}(t),
\end{align*}
where $e_{n}(t)$ is the characteristic function for the time subinterval and $\zeta_{n}$ is defined as
\begin{align*}
\zeta_{n}=\frac{1}{\sqrt{\tau}}\int_{t_{n-1}}^{t_{n}}d\zeta(t)=\frac{1}{\sqrt{\tau}}(\zeta(t_{n})-\zeta(t_{n-1}))\sim N(0,1).
\end{align*}

It should be emphasized that the exact solution of this problem is not known explicitly. We replace the ``true'' solution $u(t_{n})$ by a numerical approximation computed by small time step size $\tau^{*}=2^{-12}$ and spatial mesh size $h^{*}=1/500$, so that the temporal (or spatial) discretization error is negligible. We measure the error $e^{n}(\tau,h):=u_{h}^{n}-u(t_{n})$ by the normalized error $\mathbf{E}\|e^{n}\|^{2}$, where the expected values $\mathbf{E}\|\cdot\|$ are calculated by the averages over 100 samples. To examine the spatial and temporal convergence order separately, the convergence rates in time and space in the sense of the $L_{2}$-norm are defined as:
\begin{align*}
\mathrm{Rate}=\begin{cases}
\frac{\mathrm{ln}(\mathbf{E}\|e^{n}(\tau_{1},h)\|^{2}/\mathbf{E}\|e^{n}(\tau_{2},h)\|^{2})}{\mathrm{ln}(\tau_{1}/\tau_{2})},~\mathrm{in ~ time},\\
\frac{\mathrm{ln}(\mathbf{E}\|e^{n}(\tau,h_{1})\|^{2}/\mathbf{E}\|e^{n}(\tau,h_{2})\|^{2})}{\mathrm{ln}(h_{1}/h_{2})},~\mathrm{in ~ space},\\
\end{cases}
\end{align*}
where $\tau,\tau_{1},\tau_{2}(\neq \tau_{1})$ and $h,h_{1},h_{2}(\neq h_{1})$ are the time and space step sizes, respectively.

To begin with, all the numerical results are evaluated at $T=1$ in the given Tables. In Table 1, the numerical errors and convergence rates in spatial direction are shown by taking different spatial mesh sizes, where the fixed and sufficiently small time step sizes are taken.  It is clear that the optimal order of error estimates in space are closer to $\mathcal{O}(h^{4})$, which is in agreement with the theoretical results. Table 2 lists the numerical errors and convergence rates in temporal direction with different $\alpha$, where the spatial mesh size is fixed to be sufficiently small to ensure that dominant numerical errors come form the approximation of time-fractional derivatives. One can note that the rates of convergence are closer to the theoretical convergence order $\mathcal{O}(\tau^{2\alpha})$, that is, the numerical results are in consistent with the theoretical results.

\begin{table}[h]
\centering
\caption{Numerical errors and convergence rates in spatial direction with $\tau=2^{-10}$.}
\begin{tabular}{l l l l l l l }
\hline
\hline
\multirow{2}{*}{$h$} &
\multicolumn{2}{l}{$\alpha=0.55$}&
\multicolumn{2}{l}{$\alpha=0.70$}&
\multicolumn{2}{l}{$\alpha=0.90$}\\
\cline{2-3}
\cline{4-5}
\cline{6-7}
& $\mathbf{E}\|e^{n}\|^{2}$ & Rate &  $\mathbf{E}\|e^{n}\|^{2}$ & Rate
& $\mathbf{E}\|e^{n}\|^{2}$ & Rate\\
\hline
$1/5$  &3.4531e-02 & -     &3.4135e-02 & -    & 3.3997e-02 & -      \\
$1/10$ &2.3785e-03 & 3.86  &2.3350e-03 & 3.87 & 2.3096e-03 & 3.88 \\
$1/20$ &1.5930e-04 & 3.88  &1.5530e-04 & 3.89 & 1.5467e-04 & 3.89 \\
$1/40$ &1.1049e-05 & 3.87  &1.0477e-05 & 3.89 & 1.0010e-05 & 3.91 \\
\hline
\hline
\end{tabular}
\end{table}

\begin{table}[h]
\centering
\caption{Numerical errors and convergence rates in temporal direction with $h=1/200$.}
\begin{tabular}{l l l l l l l }
\hline
\hline
\multirow{2}{*}{$\tau$} &
\multicolumn{2}{l}{$\alpha=0.55$}&
\multicolumn{2}{l}{$\alpha=0.70$}&
\multicolumn{2}{l}{$\alpha=0.90$}\\
\cline{2-3}
\cline{4-5}
\cline{6-7}
& $\mathbf{E}\|e^{n}\|^{2}$ & Rate &  $\mathbf{E}\|e^{n}\|^{2}$ & Rate
& $\mathbf{E}\|e^{n}\|^{2}$ & Rate\\
\hline
$2^{-4}$ &2.3268e-02 & -     &2.2850e-02 & -    & 2.4233e-02 & -      \\
$2^{-5}$ &1.1379e-02 & 1.12  &8.6596e-03 & 1.40 & 6.9117e-03 & 1.81 \\
$2^{-6}$ &4.8577e-03 & 1.13  &3.2357e-03 & 1.41 & 1.9171e-03 & 1.83 \\
$2^{-7}$ &2.2197e-03 & 1.13  &1.1927e-03 & 1.42 & 5.2815e-04 & 1.84 \\
\hline
\hline
\end{tabular}
\end{table}

\section{Conclusions and discussions}

In this paper, the regularity properties of mild solution to time-fractional stochastic heat equation driven by multiplicative noise are discussed and proved. The semi-discrete and fully discrete finite element methods are developed for solving this time-fractional SPDEs. The schemes employ a standard Galerkin finite element method in space and the time direction is approximated by the GMMP scheme. The convergence error estimates for both semi-discrete and fully discrete schemes in $L_{2}$-norm are obtained. We present the numerical experiment to illustrate the accuracy of schemes, and the result fully verify the convergence theory. Noted that we only consider the Dirichlet boundary condition in our given problem, the singular boundary method might be used to deal with the complex boundary condition in our future study, because the merits of this method only needs to place the source points on the real physical boundary and uses the fundamental solutions as the kernel function [37,38,39]. In addition, numerical investigations on irregular domain problems to test the methodology potential for more realistic situations are  interesting direction for our future research.

\section*{Acknowledgements}

We would like to thank the reviewers for giving us constructive comments and suggestions
which would help us to improve the quality of the paper. This work is supported by National Nature Science Foundation of China (Grant No.11626085).

\newpage

\section*{References}

[1] H.M. Srivastava, J.J. Trujillo, Theory and applications of fractional differential equations, Elsevier, 2006.

[2] Y.Z. Povstenko, Fractional heat conduction equation and associated thermal stress, J. Therm. Stresses 28 (2004) 83-102.

[3] Z.Q. Chen, K.H. Kim, P. Kim, Fractional time stochastic partial differential equations, Stoch. Process. Appl. 125 (2015) 1470-1499.

[4] J.B. Mijena, E. Nane, Space-time fractional stochastic partial differential equations, Stoch. Proc. Appl. 125 (2015) 3301-3326.

[5] J.B. Mijena, E. Nane, Intermittence and space-time fractional stochastic partial differential equations, Potential Anal. 44 (2016) 295-312.

[6] M. Foondun, E. Nane, Asymptotic properties of some space-time fractional stochastic equations, Math. Z. (2015) 1-27.

[7] L. Chen, G. Hu, Y. Hu, J. Huang, Space-time fractional diffusions in Gaussian noisy environment, Stochastics (2016) 1-36.

[8] C. Pr\'{e}v\^{o}t, M. R\"{o}ckner, A concise course on stochastic partial differential equations, Springer, 2007.

[9] M.M. El-Borai, Some probability densities and fundamental solutions of fractional evolution equations, Chaos Soliton. Fract. 14 (2002) 433-440.

[10] J.R. Wang, Y. Zhou, Existence and controllability results for fractional semilinear differential inclusions, Nonlinear Anal. Real World
Appl. 12 (2011) 3642-3653.

[11] R. Sakthivel, S. Suganya, S.M. Anthoni, Approximate controllability of fractional stochastic evolution equations, Comput. Math. Appl. 63 (2012) 660-668.

[12] Y. Zhou, F. Jiao, Existence of mild solutions for fractional neutral evolution equations, Comput. Math. Appl. 59 (2010) 1063-1077.

[13] V. Thom\'{e}e, Galerkin finite element methods for parabolic problems, Springer-Verlag, Berlin, 1984.

[14] R. Kruse, Strong and weak approximation of semilinear stochastic evolution equations, Springer, 2014.

[15] R. Gorenflo, F. Mainardi, D. Moretti, P. Paradisi, Time fractional diffusion: a discrete random walk approach, Nonlinear Dynam. 29 (2002) 129-143.

[16] F. Zeng, C. Li, F. Liu, I. Turner, The use of finite difference/element approaches for solving the time-fractional subdiffusion equation, SIAM J. Sci. Comput. 35 (2013) A2976-A3000.

[17] L. Galeone, R. Garrappa, Explicit methods for fractional differential equations and their stability properties, J. Comput. Appl. Math. 228 (2009) 548-560.

[18] N. Ford, J. Xiao, Y. Yan, A finite element method for time fractional partial differential equations, Fract. Calc. Appl. Anal. 14 (2011) 454-474.

[19] B. Jin, R. Lazarov, Z. Zhou, Error estimates for a semidiscrete finite element method for fractional order parabolic equations, SIAM J. Numer. Anal. 51 (2013) 445-466.

[20] F. Zeng, C. Li, F. Liu, I. Turner, Numerical algorithms for time-fractional subdiffusion equation with second-order accuracy, SIAM J. Sci. Comput. 37 (2015) A55-A78.

[21] S.A. Atallah, A finite element method for time fractional partial differential equations, University of Chester, United Kingdom, 2011.

[22] B. Jin, R. Lazarov, J. Pasciak, Z. Zhou, Error analysis of semidiscrete finite element methods for inhomogeneous time-fractional diffusion, IMA J. Numer. Anal. 35 (2015) 561-582.

[23] G. Zou, B. Wang, Stochastic Burgers' equation with fractional derivative driven by multiplicative noise, Comput. Math. Appl. 74 (2017) 3195-3208.

[24] G. Zou, B. Wang, Y. Zhou, Existence and regularity of mild solutions to fractional stochastic evolution equations, Math. Model. Nat.
Phenom. (2018) https: //doi.org/10.1051/
mmnp/2018004.

[25] P.M. De Carvalho-Neto, P. Gabriela, Mild solutions to the time fractional Navier-Stokes equations in $R^{N}$, J. Differential Equations 259 (2015) 2948-2980.

[26] F. Zeng, C. Li, F. Liu, I. Turner, The use of finite difference/element approaches for solving the time-fractional subdiffusion equation, SIAM J. Sci. Comput. 35 (2013) A2976-A3000.

[27] M. Stynes, E. O'Riordan, J. L. Gracia, Error analysis of a finite difference method on graded meshes for a time-fractional diffusion equation, SIAM J. Numer. Anal. 55(2) (2017) 1057-1079.

[28] M. Stynes, J. L. Gracia, Preprocessing schemes for fractional-derivative problems to improve their convergence rates, Appl. Math. Lett. 74 (2017) 187-192.

[29] H. Wang, N. Du, A fast finite difference method for three-dimensional time-dependent space-fractional diffusion equations and its efficient implementation, J. Comput. Phys.
253(15) (2013) 50-63.

[30] F. Zeng, C. Li, F. Liu, I. Turner, Numerical algorithms for time-fractional subdiffusion equation with second-order accuracy, SIAM J. Sci. Comput. 37(1) (2015) A55-A78.

[31] M. Zheng, F. Liu, Q. Liu, K. Burrage, M.J. Simpson, Numerical solution of the time fractional reaction-diffusion equation with a moving boundary, J Comput. Phys. 338 (2017) 493-510.

[32] S. Zhai, X. Feng, Y. He, An unconditionally stable compact ADI method for three dimensional
time-fractional convection-diffusion equation, J. Comput. Phys. 269(15) (2014) 138-155.

[33] W. Deng, Finite element method for the space and time fractional Fokker-Planck equation, SIAM J. Numer. Anal. 47(1) (2008) 204-226.

[34] W. Deng, Numerical algorithm for the time fractional Fokker-Planck equation, J. Comput. Phys. 227(2) (2007) 1510-1522.

[35] W. Deng, Short memory principle and a predictor-corrector approach for fractional differential equations, J. Comput. Appl. Math. 206(1) (2007) 174-188.

[36] Y. Li, Y. Wang, W. Deng, Galerkin finite element approximations for stochastic space-time fractional wave equations, SIAM J. Numer. Anal. 55(6) (2017) 3173-3202.

[37] J. Li, W. Chen, A modified singular boundary method for three-dimensional high frequency acoustic wave problems, Appl. Math. Model. 54 (2018) 189-201.

[38] J. Li, W. Chen, Z Fu, L Sun, Explicit empirical formula evaluating original intensity factors of singular boundary method for potential and Helmholtz problems, Eng. Anal. Bound. Elem. 73 (2016) 161-169.

[39] J. Li, W. Chen, Y Gu, Error bounds of singular boundary method for potential problems, Numer. Meth. Part. D. E. 33 (2017) 1987-2004.

[40] W. Chen, G. Pang, A new definition of fractional Laplacian with application to modeling three-dimensional nonlocal heat conduction, J. Comput. Phys. 309 (2016) 350-367.

[41] W. Chen, Y. Liang, S. Hu, H. Sun, Fractional derivative anomalous diffusion equation modeling prime number distribution, Fract. Calc. Appl. Anal. 18 (2015) 789-798.

[42] G. Zou, A. Atangana, Y. Zhou, Error estimates of a semidiscrete finite element method for fractional stochastic diffusion-wave equations, Numer. Meth. Part. D. E. (2018) DOI: 10.1002/num.22252.

[43] G. Zou, G. Lv, J. Wu, Stochastic Navier-Stokes equations with Caputo derivative driven by fractional noises, J. Math. Anal. Appl. 461(1) (2018), 595-609.

[44] C.M. Elliott, S. Larsson, Error estimates with smooth and nonsmooth data for a finite element method for the Cahn-Hilliard equation, Math. Comp. 58(198) (1992) 603-630.

\bibliography{mybibfile}

\end{document}